\documentclass[12pt,letterpaper]{article}
\usepackage{amsfonts,latexsym,amssymb,amsmath,color,epsfig,mathdots,rotating}
\newtheorem{fta}{The Fundamental Theorem of Algebra.} 
 
\newtheorem{evt}{The Extreme Value Theorem.} 
 
\newtheorem{ift}{The Implicit Function Theorem.} 
 
\newtheorem{dfn}{Definition}[section]

\newtheorem{prop}[dfn]{Proposition}
 
\newtheorem{lemma}[dfn]{Lemma}

\newcommand{\thth}{^{\text{\underline{th}}}}

\newcommand{\syl}{\mathrm{Syl}}
\newcommand{\res}{\mathrm{Res}}

\newcommand{\eps}{\varepsilon}

\newcommand{\R}{\mathbb{R}}
\newcommand{\C}{\mathbb{C}}
\newcommand{\Pro}{\mathbb{P}}
\newcommand{\N}{\mathbb{N}}

\newcommand{\cI}{\mathcal{I}}
\newcommand{\cJ}{\mathcal{J}}
\newcommand{\cL}{\mathcal{L}}
\newcommand{\cX}{\mathcal{X}}

\newcommand{\Z}{\mathbb{Z}}

\newcommand{\qed}{$\blacksquare$}

\begin{document}
\title{\mbox{}\\
\vspace{-1.6in}On the BCSS Proof of the Fundamental Theorem of Algebra}  

\author{J.\ Maurice Rojas$^1$ } 

\date{\today} 

\maketitle
\footnotetext[1]{ Department of Mathematics,  
Texas A\&M University
TAMU 3368, 
College Station, Texas 77843-3368,  
USA. 
e-mail: {\tt 
rojas@math.tamu.edu} \ ,    
Web Page: {\tt 
www.math.tamu.edu/\~{}rojas} \ .
Partially supported by NSF grants DMS-1460766 and CCF-1409020, and a 
von Neumann visiting professorship. }  

\vspace{-.3in} 
\begin{abstract} 
Section 10.4 of the 1998 Springer-Verlag book {\em Complexity and Real 
Computation}, by 
Blum, Cucker, Shub, and Smale, contains a particularly elegant proof of the 
Fundamental Theorem of Algebra: The central idea of the proof naturally leads 
to a homotopy continuation algorithm for finding the roots of univariate 
polynomials, and extends naturally to a proof of B\'{e}zout's Theorem (on the 
number of roots of systems of $n$ equations in $n$ unknowns). We present a 
more detailed version of the BCSS Proof which is hopefully useful for 
students and researchers not familiar with algebraic geometry. So while there 
are no new results in this paper, the exposition is arguably more elementary. 
Any errors here are solely the responsibility of the current author.  
\end{abstract} 

% \begin{quote} 
% {\bf Note:} {\em Lecture-wise, these notes are a super-set of 
% what Alperen Erg\"ur and I spoke on during lecture 
% on Friday, June 3, 2016.}  
% \end{quote} 

\section{Preliminaries} 
The main theorem we prove is the following: 
\begin{fta} 
A (univariate) polynomial $f\!\in\!\C[x]$ of degree $d\!\geq\!1$ 
has exactly $d$ roots in $\C$, counting multiplicities. 
\end{fta} 

\noindent 
The Fundamental Theorem of Algebra (FTA) is arguably more a 
statement from complex analysis than from algebra. The first statement of 
the FTA dates back to work of Girard in 1629 and the first correct approach 
to its proof dates back to work of D'Alembert in 1746. Perhaps the first 
proof which would be considered rigorous by today's standards is 
due to Motzkin and Ostrowski, from 1933. See \cite[Sec.\ 1.8, Pg.\ 61]{rs} 
for valuable historical notes on the FTA. 

While one can certainly prove the FTA in one line 
via Liouville's Theorem from complex analysis, the proof we present below 
(from \cite[Sec.\ 10.4]{bcss}) will have two important advantages: (1) It leads 
naturally to a homotopy continuation algorithm for finding the roots of 
univariate polynomials, and (2) the proof extends naturally to a proof 
B\'{e}zout's Theorem, on the number of roots of $n$ polynomials in $n$ 
variables. 

It will be useful to recall two equivalent definitions of 
multiplicity before we eventually introduce a third, equivalent definition: 
\begin{prop}
\label{prop:mult} 
Suppose $f\!\in\!\C[x]$ is a univariate polynomial of degree $d\!\geq\!1$ 
and $m\!\in\!\N$. Then the following two conditions on a root $\zeta\!\in\!\C$ 
of $f$ are equivalent, and are referred to as {\em $f$ having a root 
of multiplicity $m$}: \\ 
(a) $f(x)\!=\!(x-\zeta)^mg(x)$ identically, for some $g\!\in\!\C[x]$ 
not divisible by $x-\zeta$.\\ 
(b) $f(\zeta)\!=\!f'(\zeta)\!=\cdots=f^{(m-1)}(\zeta)\!=\!0\!\neq\!f^{(m)}
(\zeta)$.\\
\scalebox{.96}[1]{We also call roots of multiplicity $1$ {\em non-degenerate}, 
and roots of multiplicity $\geq\!2$ {\em degenerate}. \qed}  
\end{prop} 

\noindent 
The proof of Proposition \ref{prop:mult} is straightforward (at least 
over the complex numbers) via Taylor series.  
It is worth noting that Definition (a) extends naturally to any field, while 
Definition (b) extends naturally to any analytic function on $\C$. 

Let us now fix some notation: If $f(x)
\!=\!a_0+\cdots+a_dx^d$ and $g(x)\!=\!
b_0+\cdots+b_ex^e$ are polynomials in $\C[x]$ then we define 
{\em the Sylvester matrix of $(f,g)$ (of format $(d,e)$)} to be:  \\ 
\mbox{}\hfill $\syl_{d,e}(f,g):= 
\begin{bmatrix}
a_0 & \cdots & a_{d-1} & a_d  & 0    & \cdots & \cdots & 0 \\
   & & & \ddots & & & \ddots  & \\
0   & \cdots & \cdots & 0 & a_0 & \cdots & a_{d-1} & a_d \\
b_0 & \cdots & b_e  & 0 & 0   & \cdots & \cdots & 0 \\
  & & & \ddots &  & & \ddots & \\
0   & \cdots &   & \cdots & 0 & b_0 & \cdots & b_e  
\end{bmatrix}
\begin{matrix}
\\
\left. \rule{0cm}{.9cm}\right\}
e \text{ many shifting rows}\\
\left. \rule{0cm}{.9cm}\right\}
d \text{ many shifting rows} \\
\\
\end{matrix}$\hfill\mbox{} \\ 
We also define the {\em resultant (of format $(d,e)$) of $(f,g)$} to 
be the determinant of $\syl_{d,e}(f,g)$. 
Note that the resultant (of format $(d,d-1)$) of $(f,f')$  
factors as $a_d\Delta_d(a_0,\ldots,a_d)$ for 
some polynomial $\Delta_d\!\in\!\Z[a_0,\ldots,a_d]$, thanks to the 
last column of the underlying Sylvester matrix. 
We call $\Delta_d$ the {\em degree $d$ (univariate) discriminant polynomial} 
and study it more later. That $\Delta_d$ is in fact not identically zero is a 
consequence of the following formula. 
\begin{prop}  
\label{prop:cyclo} 
For all $d\!\geq\!1$ we have $\res_{d,d-1}(x^d-1,dx^{d-1})\!=\!(-1)^{d-1}d^d$. 
\end{prop} 

\noindent 
{\bf Proof:} One need only observe that $\syl_{d,d-1}(x^d-1,dx^{d-1})$ is 
upper-triangular, with diagonal entries $-1,\ldots,-1$ ($d-1$ times) 
and $d,\ldots,d$ ($d$ times). \qed 

\medskip 
Let us define\\  
\mbox{}\hfill 
$\displaystyle{\Sigma_d:=\left\{\left. 
(\gamma_0,\ldots,\gamma_d)\!\in\!\C^{d+1}\; 
\right| \; \gamma_d\Delta_d(\gamma_0,\ldots,\gamma_d)\!=\!0 \right\}} 
\subset \C^{d+1}$.  
\hfill \mbox{} \\ 
The set $\Sigma_d$ is a slight variant of the classical discriminant variety 
$\nabla_{\{0,\ldots,d\}}\subset\Pro^d_\C$, which consists of degree $d$ 
bivariate homogeneous polynomials having a degenerate root. In particular, 
it will be convenient to identify the polynomial 
$\gamma_0+\cdots+\gamma_dx^d\!\in\!\C[x]$ with the 
point $(\gamma_0,\ldots,\gamma_d)\!\in\!\C^{d+1}$. In the Appendix, we will 
prove the following useful algebraic result: 
\begin{lemma} 
\label{lemma:disc} 
Following the notation above, assume $d,e\!\geq\!1$ and 
$\res_{d,e}(f,g)\!\neq\!0$. Then the \mbox{$2\times 1$}  
system of equations $f(x)\!=\!g(x)\!=\!0$ has {\em no} complex solutions. 
In particular, if\linebreak 
$\Delta_d(\gamma_0,\ldots,\gamma_d)\!\neq\!0$, then any complex root 
of $\gamma_0+\cdots+\gamma_dx^d$ has multiplicity $1$.
\end{lemma} 

It will also be useful to 
recall the following basic fact, easily provable via polynomial division. 
\begin{prop}
\label{prop:weak}  
\scalebox{.98}[1]{If $f\!\in\!\C[x]$ has degree $d\!\geq\!1$ then $f$ has 
at most $d$ distinct complex roots. \qed} 
\end{prop} 

\noindent 
Proposition \ref{prop:weak} (or rather, the obvious extension where $\C$ is 
replaced by an arbitrary field $F$) is sometimes called the {\em Weak} 
Fundamental Theorem of Algebra. 

We will also need some results on building paths avoiding complex hypersurfaces 
and bounding roots of polynomials: 
\begin{prop} 
\label{prop:con} 
If $g\!\in\!\C[x_1,\ldots,x_N]\setminus\{0\}$ has complex zero set $Z$ 
then $\C^N\setminus Z$ is analytically path-connected: For any 
$p,q\!\in\!\C^N\setminus Z$ there is an analytic 
$\phi : [0,1] \longrightarrow \C^N\setminus Z$ with $\phi(0)\!=\!p$ 
and $\phi(1)\!=\!q$. 
\end{prop} 

\noindent 
{\bf Proof:} We proceed by induction. Our proof will frequently use the 
obvious identification between $\R^2$ and $\C$ via real and imaginary parts. 

By an invertible affine transformation $T : \R^2\longrightarrow \R^2$, 
the $N\!=\!1$ case can clearly be reduced to the special case 
$(p,q)\!=\!(0,1)$. So then let $\R_+$ denote the positive 
ray and consider the family of parabolae $\{H_c\}_{c\in \R_+}$ defined by 
$H_c\!:=\!\{(x,y)\!\in\!\R^2\; | \; y\!=\!cx(1-x)\}$. Note that each $H_c$ 
intersects $\{0,1\}$, and can clearly be parametrized via 
$\phi(t)\!:=\!(t,ct(1-t))$. The map $\phi$ is clearly analytic and 
satisfies $\phi(0)\!=\!(0,0)$ and $\phi(1)\!=\!(1,0)$ (and the last 
two points in $\R^2$ can clearly also be viewed as the points $0$ and $1$ 
in $\C$). 

Clearly, there are only finite many values of $c$ with 
$H_c\cap Z\!\neq\!\emptyset$: If $Z\!=\!\{(x_1,y_1),\ldots,(x_\ell,y_\ell)\}$ 
then we simply need to pick a $c$ satisfying $cx_i(1-x_i)\!\neq\!y_i$ for all 
$i$. Picking $c\!=\!\frac{1}{2}\min_i \frac{y_i}{x_i(1-x_i)}$ clearly suffices. 
We thus obtain an $H_c$ connecting $0$ and $1$, and the resulting parametrized 
path $\phi : [0,1]\longrightarrow \R^2$ is clearly analytic, as is also the 
composition $T^{-1}(\phi(\cdot))$. So we are done. 

To prove the general case we will let $N\!\geq\!2$ and reduce to the 
case $N\!=\!1$ by restricting to a cleverly embedded copy of $\C$ in 
$\C^N$: Define $L : \C\longrightarrow \C^N$ by $L(t)\!=\!(1-t)p+tq$, 
let $h(t)\!:=\!g(L(t))$, and let $W$ denote the zero set of $h$ in $\C$. 
Note that $L(0)\!=\!p$ and $L(1)\!=\!q$ and, 
by assumption, $h(0)h(1)\!=\!g(p)g(q)\!\neq\!0$. So, since $L$ is an 
analytic function, it suffices to find an analytic path in $\C\setminus W$ 
connecting $0$ and $1$. By the $N\!=\!1$ case we have already proved, we 
are done. \qed  

\begin{lemma} 
\label{lemma:bound} 
If $f(x)\!=\!\gamma_0+\cdots+\gamma_dx^d\!\in\!\C[x]$ is a univariate 
polynomial of degree $d\!\geq\!1$ then all the complex roots of $f$ lie in 
the open disk of radius $2\max\limits_{i\in\{0,\ldots,d-1\}} 
\left|\frac{\gamma_i}{\gamma_d}\right|^{\frac{1}{d-i}}$ centered at $0$. 
\end{lemma} 

\noindent 
{\bf Proof:} The proof boils down to a judicious use of the Triangle 
Inequality. First note that if $\zeta\!\in\!\C$ is a root of $f$ then 
$|\gamma_d\zeta^d|\!=\!|\gamma_0+\cdots+\gamma_{d-1}\zeta^{d-1}|$ and 
thus 
\begin{eqnarray}
\label{ineq:central} 
|\gamma_d||\zeta|^d & \leq & |\gamma_0|+|\gamma_1||\zeta|+\cdots+|\gamma_{d-1}||\zeta|^{d-1}
\end{eqnarray}  

So let us now consider the special case where $\gamma_d\!=\!1$ and 
$\max\limits_{i\in\{0,\ldots,d-1\}} |\gamma_i|\!=\!1$. In this case,  
Inequality (\ref{ineq:central}) implies that 
$|\zeta|^d \! \leq \! 1+|\zeta|+\cdots+|\zeta|^{d-1}$. So by  
geometric series we obtain $|\zeta|^d \! \leq \! \frac{|\zeta|^d-1}
{|\zeta|-1}$. In the event that $|\zeta|\!\geq\!2$ we would then obtain 
$|\zeta|^d \! \leq \! |\zeta|^d-1$, which is clearly impossible. So we 
see that $|\zeta|\!<\!2$ and our lemma is true. 

To obtain the general case, let $C\!:=\!\max\limits_{i\in\{0,\ldots,d-1\}} 
\left|\frac{\gamma_i}{\gamma_d}\right|^{\frac{1}{d-i}}$ and 
consider the polynomial 
\[ g(x)\!:=\!\frac{f(Cx)}{\gamma_d C^d}.\]  
(Note that 
$\gamma_d\!\neq\!0$ since $f$ is assumed to have degree $d$.) Clearly, 
our lemma is proved if we can prove that all the complex roots of $g$ have 
norm less than $2$. Toward this end, observe that the coefficient of 
$x^i$ in $g$ (for $i\!\in\!\{0,\ldots,d-1\}$) 
is $\frac{\gamma_iC^i}{\gamma_d C^d}\!=\!\left(
\frac{(\gamma_i/\gamma_d)^{1/(d-i)}}{C}\right)^{d-i}$. Since 
$|\gamma_i/\gamma_d|^{1/(d-i)}\!\leq\!C$ by the definition of $C$, 
every coefficient of $g$ has norm at most $1$. Furthermore, the coefficient 
of $x^d$ in $g$ is $1$. So by our early special case, we are done. \qed 

\noindent 
Lemma \ref{lemma:bound} is in fact a very special case of more general 
estimates on complex roots of univariate polynomials. Such bounds 
date back to work of Cauchy around 1829, if not earlier. 
More recent work relates the norms of complex roots to Archimedean 
tropical varieties (e.g., \cite[Thm.\ 1.7]{aknr}) and extends such 
estimates to polynomials in multiple variables. 
% In particular, the radius above is simply a 
% disguised form of $e^s$ where $s$ is the maximal slope of a lower  
% edge of $\mathrm{ArchNewt}(f)$. 

% \begin{hb} \cite{borel} 
% Any open cover $\{U_\alpha\}_\alpha$ of a closed and bounded subset 
% of $\R^N$ admits a finite subcover. \qed 
% \end{hb} 
Finally, we'll need the following two analytic results. 
\begin{evt}
(See, e.g., \cite{munkres}.) If $S\!\subset\!\R^N$ is closed and bounded and 
$f : S\longrightarrow \R$ is continuous, then there exist $a_1,a_2\!\in\!S$ 
with $f(a_1)\!=\!\inf_{x\in S} f(x)$ and $f(a_2)\!=\!\sup_{x\in S} f(x)$. In 
particular, $f$ is bounded. \qed 
\end{evt} 
\begin{ift} 
(See, e.g., \cite{hubbard}.)  
If $G : \C^2 \longrightarrow \C$ is an analytic function with 
$G(t_0,x_0)\!=\!0$ 
and $\frac{\partial G}{\partial x}|_{(t_0,x_0)}\!\neq\!0$ for some 
$(t_0,x_0)\!\in\!\C^2$, then there is an open disk $U\!\subset\!\C$ containing 
$t_0$, and an analytic function $\cX : U \longrightarrow \C$ (called a 
{\em branch of the implicit function defined by $G(t,x)\!=\!0$}) with 
$\cX(t_0)\!=\!x_0$ and $G(t,\cX(t))\!=\!0$ for all $t\!\in\!U$. \qed  
\end{ift} 

\noindent 
The Extreme Value Theorem is covered in most any basic course in real 
analysis. We also note that \cite[Sec.\ 8.2]{bcss} presents an elegant proof of 
the {\em Inverse} Function Theorem (via clever infinite series manipulations) 
that can be used to derive explicit estimates for the size of the disk 
in the Implicit Function Theorem. 

\section{The Proof of the Fundamental Theorem of Algebra} 
We first observe that the case $d\!=\!1$ is trivial: $\gamma_0+\gamma_1x$ has 
only the root $-\gamma_0/\gamma_1$, which has multiplicity $1$ since 
$f'(x)$ is identically $\gamma_1$, and $\gamma_1\!\neq\!0$ since we assumed the 
degree of $f$ is $d\!=\!1$. So we may assume $d\!\geq\!2$ henceforth. 

We then separate our proof into two cases: (1) $f\not\in \Sigma_d$ and 
(2) $f\!\in\!\Sigma_d$. The first case is harder so we take care of it 
first. 

% \begin{quote}
% {\bf Note:} {\em In what follows, since we are identifying polynomials 
% of degree $\leq\!d$ with points in $\C^{d+1}$, we will be speaking of 
% polynomials as points, and vice-versa, wherever convenient.} 
% \end{quote} 

\subsection{The Case $f\not \in \Sigma_d$} First observe that 
$x^d-1\not\in\Sigma_d$ thanks to Proposition \ref{prop:cyclo}.
Since $\Sigma_d$ is the complex zero set of a polynomial,                
Proposition \ref{prop:con} tells us that
there is an analytic function\\ 
\mbox{}\hfill $\phi\!=\!(\phi_0,\ldots,\phi_d) : 
[0,1]\longrightarrow \C^{d+1}\setminus
\Sigma_d$\hfill\mbox{}\\ 
with $\phi(0)\!=\!x^d-1$ and $\phi(1)\!=\!f$. 
(Since we identify polynomials and points, we could have also said 
$\phi(0)\!=\!(-1,0,\ldots,0,1)$ and $\phi(1)\!=\!(\gamma_0,\ldots,\gamma_d)$.) 
Let $G : [0,1]\times \C \longrightarrow \C$ be the analytic function 
defined by $G(t,x)\!:=\!\phi_0(t)+
\cdots+\phi_d(t)x^d$, i.e., $G$ simply evaluates the polynomial $\phi(t)$ 
on our path at the complex number $x$. 

Our main trick then will be to use our path of polynomials 
to build $d$ {\em paths of roots}. On one end ($t\!=\!0$), the root paths 
start at the $d\thth$ roots of unity. We will ultimately prove that, on the 
other end ($t\!=\!1$), our root paths yield $d$ distinct complex roots for $f$. 
Since $f\not\in\Sigma_d$ implies that each 
root of $f$ has multiplicity $1$, and Proposition \ref{prop:weak} tells us 
that $f$ can have no more than $d$ complex roots, proving that $f$ 
has $d$ complex roots will complete our proof.  

Letting $D\subseteq[0,1]$ 
denote the set of $t\!\in\![0,1]$ such that $Z_t\!:=\!\{\zeta\!\in\!\C\; | \; 
G(t,\zeta)\!=\!0\}$ has cardinality $d$, it clearly suffices to prove 
$D\!=\![0,1]$. Toward this end, we'll first establish the following facts: \\  
\mbox{}\hfill 
($P_0$) $D$ is non-empty. \  \ \ \ \   
($P_1$) $D$ is open. \ \ \ \ \ 
($P_2$) $D$ is closed.  \hfill\mbox{} 
 
\scalebox{.94}[1]{To prove ($P_0$), first note that 
$x^d-1$ has exactly $d$ roots: They are 
$e^{2\pi\sqrt{-1}/d},\ldots,e^{2\pi(d-1)\sqrt{-1}/d}$,}\linebreak 
thanks to Euler's 1748 formula for the exponential function. The Implicit 
Function\linebreak 
Theorem applied to each of these roots then implies that 
there are a positive $\eps$, and analytic functions 
$\cX_1,\ldots,\cX_d : [0,\eps] \longrightarrow \C$, with 
$\{\cX_1(0),\ldots,\cX_d(0)\}\!=\!\left\{e^{2\pi\sqrt{-1}/d},\ldots,
e^{2\pi(d-1)\sqrt{-1}/d}\right\}$ and $G(t,\cX_i(t))\!=\!0$ for all 
$t\!\in\![0,\eps)$ and $i\!\in\!\{1,\ldots,d\}$. Furthermore, 
since the $\cX_i$ are\linebreak 
continuous and $|\cX_i(0)-\cX_j(0)|\!>\!0$ for all 
$i\!\neq\!j$, we can assume further that $|\cX_i(t)-\cX_j(t)|\!>\!0$ for all 
$i\!\neq\!j$ and $t\!\in\![0,\eps')$, for some positive $\eps'\!\leq\!\eps$. 
Since $Z_t$ has cardinality $\leq\!d$ 
for all $t$ (thanks to Proposition \ref{prop:weak}), we thus obtain 
that $D\!\supseteq\![0,\eps')$. So ($P_0$) is proved, and we refer 
to analytic functions $\cX_i$ satisfying $G(t,\cX_i(t))$ 
for $t$ in some open interval as {\em root branches}. (So the 
root paths we are trying to build are connected subsets of the graphs 
of the root branches.)  

The proof of Statement ($P_1$) follows the proof of Statement ($P_0$) almost 
identically: If $Z_t$ has 
cardinality $d$ for some $t\!\in\![0,1]$,  
then the Implicit Function Theorem applied to each of the $\zeta\!\in\!Z_t$ 
yields $d$ root branches defined in an open neighborhood $U\!\ni\!t$. By 
continuity again, we obtain that there is a non-empty open set 
$U'\!\subseteq\!U$, with $U'\!\ni\!t$ as well, such that the root paths do  
not intersect at any $t\!\in\!U'$. So then $U'\!\subseteq\!D$  
and ($P_1$) is proved. 

To prove Statement ($P_2$) we first need to observe that 
$|\phi_0|,\ldots,|\phi_d|$ are continuous\linebreak functions 
on the compact set $[0,1]$. So by the Extreme Value Theorem, there is 
a\linebreak 
positive $M$ with $|\phi_0(t)|,\ldots,|\phi_{d-1}(t)|\!\leq\!M$ for all 
$t\!\in\![0,1]$. Furthermore, 
$\phi_d$ is nonzero on all of $[0,1]$ (by the definition of $\Sigma_d$) 
and thus $|\phi_d|$ is a positive on $[0,1]$. 
So by the Extreme Value Theorem once again, $|\phi_d|$ has a positive minimum 
$m$ on $[0,1]$, and thus $\left|\frac{\phi_i(t)}{\phi_d(t)}\right|\!\leq\!M/m$ 
for all $i\!\in\!\{0,\ldots,d-1\}$ and $t\!\in\![0,1]$. 
By Lemma \ref{lemma:bound} we thus obtain that\\ 
\mbox{}\hfill  
$\cI:=\{(t,\zeta)\!\in\![0,1]\times \C \; | \; G(t,\zeta)\!=\!0\}$\hfill
\mbox{}\\  
is a bounded set. $\cI$ is also closed since it is the intersection 
of a closed (and bounded) rectangle with the zero set of a continuous 
function. 

Suppose now that $t'\!\in\![0,1]\setminus D$ is a limit point of 
$D$. Then there are sequences\\
\mbox{} \hfill $\zeta^{(1)}_1,\zeta^{(1)}_2, \ldots \in \C$, \hfill\mbox{}\\ 
\mbox{} \hfill $\vdots$ \hfill\mbox{}\\ 
\mbox{} \hfill $\zeta^{(d)}_1,\zeta^{(d)}_2, \ldots \in \C$, \hfill\mbox{}\\ 
and $t_1,t_2,t_3,\ldots\in D$,  
such that $\left\{\zeta^{(1)}_i,\ldots,\zeta^{(d)}_i\right\}$ has cardinality 
$d$ for all $i\!\in\!\N$, $G\!\left(t_i,\zeta^{(j)}_i\right)\!=\!0$ 
for all $(i,j)\!\in\!\N\times \{1,\ldots,d\}$, and 
$\lim\limits_{i\rightarrow \infty}t_i\!=\!t'$. 

We claim that 
the sequences of roots $\left(\zeta^{(1)}_i\right)_{i\in\N},\ldots,
\left(\zeta^{(d)}_i\right)_{i\in\N}$ always remain at a positive distance from 
each other. More precisely, we claim that there is a $\delta$ 
with $\left|\zeta^{(j)}_i-\zeta^{(j')}_i\right|\!\geq\!\delta\!>\!0$ for all 
$j\!\neq\!j'$ and $i\!\in\!\N$. For if not, then there is a pair $(j,j')$ 
with $j\!\neq\!j'$ and a subsequence $i_1,i_2,\ldots$ with 
$\lim_{n\rightarrow\infty}\left|\zeta^{(j)}_{i_n}-\zeta^{(j')}_{i_n}\right|
\!=\!0$, i.e., both $\zeta^{(j)}_{i_n}$ and $\zeta^{(j')}_{i_n}$ tend 
to a common limit $\zeta$. Since $G$ and $\frac{\partial G}{\partial x}$ are 
continuous, we would then obtain $G(t',\zeta)\!=\!0$ and 
\[
\left. \frac{dG}{dx}\right|_{(t',\zeta)}=
\lim_{n \rightarrow \infty}
\frac{G\left(t_{i_n},\zeta^{(j)}_{i_n}\right)-G\left(t_{i_n},
\zeta^{(j')}_{i_n}\right)}
{\zeta^{(j)}_{i_n}-\zeta^{(j')}_{i_n}} =
\lim_{n \rightarrow \infty}
\frac{0-0}{\zeta^{(j)}_{i_n}-\zeta^{(j')}_{i_n}} =
0. \] 
But then, this contradicts the fact that $\phi([0,1])$ avoids $\Sigma_d$. 
So the sequences $\left(\zeta^{(1)}_i\right)_{i\in\N},\ldots,$\linebreak  
$\left(\zeta^{(d)}_i\right)_{i\in\N}$ are indeed separated by some positive 
$\delta$.  

By this separation, and from the fact that $\cI$ is closed and bounded, 
we then obtain that each sequence $\left(\zeta^{(j)}_i\right)_{i\in \N}$ has 
a limit point $\zeta^{(j)}$, and these limit points are all separated by 
a distance at least $\delta$ from each other. In particular, again using the 
continuity of $G$, we obtain $G\!\left(t',\zeta^{(j)}\right)\!=\!0$ for all 
$j\!\in\!\{1,\ldots,d\}$. So $Z_{t'}$ has cardinality $\geq\!d$ and thus, 
by Proposition \ref{prop:weak} again, $t'\!\in\!D$. So we at last have 
Statement ($P_2$). 

To conclude, note that we must have $D\!=\![0,1]$ since $[0,1]$ is 
connected, and $D$ being non-empty, open, and closed implies that $D$ must 
be a connected component of $[0,1]$. \qed 

\subsection{The Case $f \in \Sigma_d$} First note that we may assume 
$\gamma_d\!\neq\!0$, for otherwise we could simply decrease $d$ and 
reduce to a lower degree case. 

Now let $L(t,x)\!:=\!(1-t)(\gamma_dx^d-\gamma_d)
+tf(x)$ and observe that the coefficient of $x^d$ in $L(t,x)$ is exactly 
$\gamma_d\!\neq\!0$ for all $t\!\in\!\C$. So the roots of 
$L(0,x)$ are the roots of unity, $L(1,x)\!=\!f(x)$, and $L(t,x)$ has 
degree $d$ in $x$ for all $t\!\in\![0,1]$. Furthermore, since 
$x^d-1\not\in \Sigma_d$ (as already observed in our last case), 
$\Delta_d(L(t,x))$ is a polynomial of {\em positive} degree $\leq\!2d-1$ 
in $t$. So, by restricting $L(t,\cdot)$ to $[\eps,1]$ for some suitable 
$\eps\!\geq\!0$ and rescaling $t$, we can build an analytic path 
$\psi\!=\!(\psi_0, \ldots,\psi_d) : [0,1]\longrightarrow \C^{d+1}$ 
with $\psi(0)\!=\!g$ for some $g\!\in\!\C^{d+1}\setminus\Sigma_d$ and 
{\em only} intersecting $\Sigma_d$ at $\psi(1)\!=\!f$. We will use 
$\psi$ to build root paths that prove that $f$ has at least one root in 
$\C$ and then conclude by long division. 

More precisely, similar to our last case, let \\ 
\mbox{}\hfill 
$H(t,x)\!:=\!\psi_0(t)+\cdots+\psi_d(t)\zeta^d$ 
\hfill\mbox{}\\ 
\mbox{}\hfill 
$W_t\!:=\!\{\zeta\!\in\!\C\; | \; H(t,\zeta)\!=\!0\}$ 
\hfill\mbox{}\\ 
\mbox{}\hfill 
$\cJ:=\{(t,\zeta)\!\in\![0,1]\times \C \; | \; H(t,\zeta)\!=\!0\}$. 
\hfill\mbox{}\\ 
Thanks to the framework of the $f\not\in\Sigma_d$ case (which we've 
already proved), it is clear that 
$H$ is analytic, $\cJ\setminus (\{1\}\times \C)$ is a disjoint union of 
$d$ real analytic curves, and $\cJ$ is closed and bounded. In particular,  
these curves induce (possibly overlapping) limit points 
$(1,\rho_1),\ldots,(1,\rho_d)$. Furthermore, since $H$ is continuous, 
$f$ must vanish at each $\rho_i$. So $f$ has at least one root $\rho\!\in\!\C$ 
and, letting $\mu$ ($\geq\!1$) denote the multiplicity of $\rho$, 
it suffices to prove that $\frac{f(x)}{(x-\rho)^\mu}$ has exactly 
$d-\mu$ roots in $\C$ counting multiplicities. So by induction, we are done. 
\qed 

\section{Consequences for Intersection Multiplicity and\\ Resultants} 
We first point out that the preceding proof yields a {\em geometric} 
definition of intersection 
multiplicity: For a degenerate root $\rho\!\in\!\C$ of $f$, 
following the notation of the last section, 
we can define the intersection multiplicity of a root $\rho$ to be 
the number of paths of $\cJ$ ending at $(1,\rho)$, i.e., the number of 
$\rho_i$ equal to $\rho$. With just a little more work (see, e.g., 
\cite[Ch.\ 10]{bcss}), we can show that this definition makes sense, 
i.e., our definition is independent of how we build a path of polynomials 
in $\C^{d+1}\setminus\Sigma_d$ tending to $f$. It turns out that this 
definition of intersection multiplicity generalizes easily to 
higher-dimensions, e.g., B\'ezout's Theorem (see \cite[Sec.\ 10.5]{bcss} 
for the details). 

A useful consequence of the Fundamental Theorem of Algebra is the 
following inverse of the implication from Lemma \ref{lemma:disc}: 
\begin{lemma} 
\label{lemma:conv} 
Following the notation of Lemma \ref{lemma:disc}, assume instead that 
$\res_{d,e}(f,g)\!=\!0$ and $a_db_e\!\neq\!0$. Then there is a 
$\zeta\!\in\!\C$ with $f(\zeta)\!=\!g(\zeta)\!=\!0$.  
\end{lemma} 

\noindent 
We prove Lemma \ref{lemma:conv} in the Appendix. 

\section{A Root Approximation Algorithm from our Proof}  
Given any $f\!\in\!\C[x]$, consider the function 
$H(t,x)\!:=\!(1-t)(x^d-1)+tf(x)$. Let $N$ be a large integer to be 
specified later, and consider the following iteration: Define 
$z^{(1)}_0,\ldots,z^{(d)}_0$ to be the $d\thth$ roots of unity, and 
then set 
\[ z^{(j)}_{n+1}:=z^{(j)}_n - \left(\left. 
                        \frac{\partial H}{\partial x}
                    \right|_{\left(\frac{n+1}{N},z^{(j)}_n\right)}\right)^{-1}
                    H\!\left(\frac{n+1}{N},z^{(j)}_n\right) \] 
for all $j\!\in\!\{1,\ldots,d\}$ and $n\!\in\!\{0,\ldots,N-1\}$. This is the 
most basic {\em homotopy continuation method} for approximating the roots of 
$f$: Under certain reasonably enforceable conditions, $z^{(1)}_N,
\ldots,z^{(d)}_N$ are high-precision approximations to all the roots of 
$f$. 

The key trick behind this method is that the sequences\\  
\mbox{}\hfill 
$\left(H\!\left(\frac{n}{N},x\right),z^{(1)}_n\right)_{n\in \N},\ldots, 
\left(H\!\left(\frac{n}{N},x\right),z^{(d)}_n\right)_{n\in \N}$\hfill\mbox{}\\ 
respectively 
yield discrete approximations to pieces of the {\em incidence variety}\\ 
\mbox{}\hfill $\cJ\!:=\!\left\{(c_0,\ldots,c_d,
\zeta)\!\in\!\C^{d+1}\times \C\; |\; c_0+\cdots+c_d\zeta^d\!=\!0\right\}$. 
\hfill\mbox{}\\   
In particular, the incidence variety consists of all pairs consisting of a 
polynomial and its roots, and we are trying to recover a point on $\cJ$ 
from its projection onto the first factor: the coefficient space $\C^{d+1}$. 
The preceding iteration takes $d$ (easily obtainable) points on 
$\cJ$ and then applies Newton's method (to deformations of $f$) to 
follow the $d$ curves leading to the (hard to find) pairs 
$(f,\zeta_1),\ldots,(f,\zeta_d)$ 
where $\left\{\zeta^{(1)},\ldots,\zeta^{(d)}\right\}$ is the multi-set of 
roots of $f$ (appearing with appropriate multiplicity). 

If by luck we have that the line segment $\cL:=\left\{ 
(1-t)(x^d-1)+tf(x)\right\}_{t\in[0,1]}$ does {\em not} intersect 
$\Sigma_d$, and $N$ is sufficiently large, then one can in fact guarantee that 
the underlying Newton iterations are all well-defined and that the 
points $z^{(1)}_N,\ldots,z^{(d)}_N$ are approximate roots of $f$ {\em 
in the sense of Smale}. 
One should of course not rely on luck. However, the assumption that 
$\cL$ does not intersect $\Sigma_d$ actually holds with probability $1$ 
for many natural probability measures on the space of coefficients 
$\C^{d+1}$. And, in the event that $\cL$ intersects $\Sigma_d$, there 
are well-known strategies to modify the underlying path and iteration. 

As for the size of the paramter $N$, there is a rigorous notion of 
{\em approximate root} 
(due to Smale \cite{smale}) that gives a practical sufficient 
condition for an approximation $z^{(j)}_N$ to be efficiently 
refinable (to arbitrarily small distance from a true root) 
via Newton's method applied directly to $\left(f,z^{(j)}_N\right)$. 
% We will make these notions precise later in lecture, but four seminal 
Four seminal 
references on analyzing how large $N$ should be are \cite{ss,bcss,bp,lairez}. 

This approach to numerical solving has a transparent generalization to 
systems of multivariate polynomial equations (see, e.g., \cite{som}), and 
has already been implemented and proven quite valuable. An excellent example 
is the package {\tt Bertini}, which is freely downloadable after a simple 
google search. Homotopy methods thus possess a certain conceptual beauty, 
and making them faster and fully reliable leads to many active research 
problems in (numerical) algebraic geometry. 

\section*{Appendix: Proving the Equivalence of Vanishing\\ Discriminants 
and the Presence of Degenerate Roots} 
Collecting powers of $x$, it is easy to prove the following 
proposition, which we'll use below: 
\begin{prop} 
\label{prop:mat} 
For any $d,e\!\geq\!1$, and complex constants $\alpha_i$, $\beta_i$, 
$\gamma_i$, and $u_i$, we have that the matrix product identity\\ 
\mbox{}\hfill $[\alpha_0,\ldots,\alpha_{e-1},\beta_0,\ldots,
\beta_{d-1}]\syl_{d,e}(f,g)\!=\![u_0,\ldots,u_{d+e-1}]$ 
\hfill \mbox{}\\  
holds if and only if the polynomial identity \\ 
\mbox{}\hfill 
$(\alpha_0+\cdots+\alpha_{e-1}x^{e-1})f+(\beta_0+\cdots+\beta_{d-1}x^{d-1})
g(x)=u_0+\cdots+u_{d+e-1}x^{d+e-1}$\hfill\mbox{} \\ 
holds. \qed  
\end{prop} 

\subsection*{The Proof of Lemma \ref{lemma:disc}} 
If $d\!=\!e\!=\!1$ then the resultant not vanishing is equivalent to 
the vectors $[a_0,a_1]$ and $[b_0,b_1]$ being linearly independent. 
So it is clear that $f\!=\!g\!=\!0$ can have no complex solutions. So 
let us assume, swapping $f$ and $g$ if necessary, that $d\!\geq\!2$. 

Then $\res_{d,e}(f,g)\!\neq\!0$ implies that $\det\syl_{d,e}(f,g)\!\neq\!0$,  
and thus the linear system\\ 
\mbox{}\hfill $[r_0,\ldots,r_{e-1},
s_0,\ldots,s_{d-1}]\syl_{d,e}(f,g)\!=\![1,0,\ldots,0]$
\hfill\mbox{}\\ 
has a unique solution. 

By Proposition \ref{prop:mat} we then have that \\  
\mbox{}\hfill $(r_0+\cdots+r_{e-1}x^{e-1})f(x)+(s_0+\cdots+s_{d-1}
x^{d-1})g(x)\!=\!1$ \hfill \mbox{}\\  
identically. This implies, in particular, that 
$f$ and $g$ can have no complex roots in common: Evaluating our last 
polynomial linear combination at such a root would imply that $0\!=\!1$ 
--- a contradiction. So we obtain the main assertion. 

The consequence for $\Delta_d$ follows immediately upon setting 
$g\!=\!f'$ and observing that any root $\zeta\!\in\!\C$ of $f$ with 
$f'(\zeta)\!\neq\!0$ is a root of multiplicity $1$, thanks to 
Proposition \ref{prop:mult}. \qed 
 
\subsection*{The Proof of Lemma \ref{lemma:conv}} 
Under the stated hypotheses, $\det \res_{d,e}(f,g)\!=\!0$. 
So Proposition \ref{prop:mat} implies that 
there must be polynomials $r,s\!\in\!\C[x]$, at least one of which 
is not identically $0$, with degrees respectively 
bounded from above by $e-1$ and $d-1$, satisfying $rf+sg\!=\!0$ 
identically. In particular, the factorizations of $rf$ and $-sg$ must be 
identical. Note also that if $s$ is identically $0$ then $f$ must be 
identically $0$, and we trivially have that $f$ vanishes at all the roots 
of $g$. (And $g$ has at least $1$ complex root, thanks to the Fundamental 
Theorem of Algebra.) So we may assume that {\em both} $r$ and $s$ are not 
identically $0$. 

By the Fundamental Theorem of Algebra, 
$f$ has exactly $d$ degree $1$ factors (counting 
multiplicity) since $a_d\!\neq\!0$. Since $s$ has no more than 
$d-1$ linear factors, $f$ and $g$ must have a degree $1$ factor in common. 
So we are done. \qed  

\bibliographystyle{amsalpha}

\end{document}